\documentclass{amsart}
\usepackage{amssymb,latexsym,graphics,amsfonts}

\swapnumbers
 \theoremstyle{plain}

 \theoremstyle{definition}
 
 \theoremstyle{remark}

 \numberwithin{equation}{subsection}

 \newcommand{\cal}[1]{\mathcal{#1}}

\newcommand{\ti}{\times}

\newcommand{\subs}{\subseteq}
\newcommand{\al}{\alpha}

\newcommand{\la}{\lambda}

\newcommand{\lo}{\longrightarrow}
\newcommand{\Lo}{\Longrightarrow}
\newcommand{\LO}{\Longleftrightarrow}

\newcommand{\BN}{\Bbb N}

\newcommand{\cA}{\ensuremath{\mathcal A}}

\begin{document}

\title{The fourth duals of  Banach algebras}

\author{M. Eshaghi Gordji}
\address{Department of Mathematics,
University of Semnan, Semnan, Iran} \email{madjideg@walla.com}

\author{S. A. R. Hosseiniun}
\address{Department of Mathematical Sciences,
Shahid Beheshti University, Tehran, Iran}
\email{ahosseinioun@yahoo.com}



\dedicatory{}



\smallskip

\begin{abstract}
Let $\cal A$ be a Banach algebra. Then $\cal A^{**}$ the second dual
of $\cal A$ is a Banach algebra with first (second) Arens product.
We study the Arens products of $\cal A^{4}(=({\cal A^{**}})^{**}).$
We found some conditions on $\cal A^{**}$ to be a left ideal in
$\cal A^{4}.$ We found the biggest two sided ideal $I$ of ${\cal
A},$ in which $I$ is a left (right) ideal of ${\cal A}^{**}$.
\end{abstract}

\maketitle


\section{Introduction}

The regularity of bilinear maps on norm spaces, was introduced by
Arens in 1951 [1]. Let $X$ , $Y$ and $Z$ be normed spaced and let
$f:X\times Y\longrightarrow Z$ be a continuous
 bilinear map, then
$f^*:Z^*\times X\longrightarrow Y^*$ (the transpose of $f$ ) is
defined by
$$\langle f^*(z^*,x),y\rangle =\langle z^*,f(x,y)\rangle
\hspace{1cm}(z^*\in Z^*~,~x\in X~,~ y\in Y).$$ ($f^*$ is a
continuous bilinear map). Clearly, for each $x\in X$, the mapping
$z^*\longmapsto f^*(z^*,x):Z^*\longrightarrow Y^*$ is
$weak^*-weak^*$ continuous. We take $f^{**}=(f^*)^*$ and
$f^{***}=(f^{**})^*$,... .

 Let $X,Y$ and $Z$ be Banach spaces and
let $f:X\ti Y\lo Z$ be a bilinear map. Let $X$ and $Z$ be dual
Banach spaces then we define
\begin{equation*}
Z_r(f):=\{y\in Y: f(.,y): X\lo Z\; \text{~~~~is }\;
weak^*-weak^*-continuous\}\hspace {1.2cm}(1.1),
\end{equation*}
so if $Y$ and $Z$ are dual spaces then we define
\begin{equation*}
Z_l(f):=\{x\in X: f(x,.): Y\lo Z \;\text{~~~~is}\;
weak^*-weak^*-continuous\}\hspace {1.2cm}(1.2).
\end{equation*}
We call $Z_r(f)$ and $Z_l(f),$ the topological centers of $f$. For
example if ${\cal A}$ is a Banach algebra by product $\pi:{\cal
A}\ti {\cal A}\lo {\cal A}$, then ${\cal A}^{**}$ the second dual of
${\cal A}$ is a Banach algebra by each of products $\pi^{***}$ and
$\pi^{r***r}$ (this products are  the first and the second Arens
products of ${\cal A}^{**}$ respectively). Also we have
$Z_l(\pi^{***})=Z_1$ the left topological center of ${\cal A}^{**}$
and $Z_r(\pi^{r***r})=Z_2$ the right topological center of ${\cal
A}^{**}$ (see [2], [4], [7]).

 Let $Z=\{a''\in {\cal A}^{**}: \pi^{***}(a'',b'')=\pi^{r***r}
(a'',b'')\; \text{for every}\; b''\in {\cal A}^{**}\}$. Then it is
easy to show that ${\cal A}\subs Z$ if and only if  ${\cal A}$ is
commutative.
\paragraph{\bf Lemma 1.1.} Let ${\cal A}$ be a commutative Banach algebra.
Then the following assertions hold.

(i) If ${\cal A}^{**}$ has identity $E$ for one of the Arens
products then $E$ is identity for other product.

(ii) For every $a''\in {\cal A}^{**},$ we have
$\pi^{***}(a'',a'')=\pi^{r***}(a'',a'')=\pi^{r***r}(a'',a'')$.
\paragraph{\bf Proof.} (i) Let $E$ be the identity for $({\cal A}^{**},\pi^{***}).$ Then for
every $F\in {\cal A}^{**}$, we have
$$\pi^{***r}(F,E)=\pi^{***}(E,F)=F=\pi^{***}(F,E)=\pi^{r***}(F,E)=\pi^{r***r}(E,F).$$
Similarly we can show that $E$ is  the identity for $({\cal
A}^{**},\pi^{***})$ when $E$ is the identity of $({\cal
A}^{**},\pi^{r***r}).$ The proof of (ii) is
trivially.\hfill$\blacksquare~$

 In this paper we study the Arens regularity of ${\cal
A}^{(4)}.$ We fined the conditions on $\cal A^{**}$ to be a left
ideal in $\cal A^{****}.$ Finally we fined the biggest two sided
ideal $I$ of ${\cal A}$ in which $I$ is a left (right) ideal of
${\cal A}^{**}$.

\section{Arens products of  ${\cal A}^{****}$ }

Let ${\cal A}$ be a Banach algebra by product $\pi:{\cal A}\ti {\cal
A}\lo {\cal A}$.  The fourth dual ${\cal A}^{(4)}=({\cal
A}^{**})^{**}$ of $\cal A,$ is a Banach algebra by products
$\pi^{******},\pi^{***r***r}, \pi^{r***r***}$ and $\pi^{r******r}$.
Also it is easy to show that
$$Z_r(\pi^{******})=Z_l(\pi^{***r***r})=Z_r(\pi^{r***r***})=Z_l(\pi^{r******r})={\cal
A}^{(4)}.$$
\paragraph{\bf Theorem 2.1.} Let $\cal A$ be a Banach algebra,
then the following assertions are equivalent

(i)  $Z_l(\pi^{r******r})={\cal A}^{(4)}$

(ii) $Z_l(\pi^{***r***r})={\cal A}^{(4)}$

(iii) $Z_r(\pi^{******})={\cal A}^{(4)}$

(iv) $Z_r(\pi^{r***r***})={\cal A}^{(4)}.$
\paragraph{\bf Proof.} We have
$$ Z_l(\pi^{******})=\cal A^{(4)}\LO (\cal A^{**},\pi^{***})\;\text{is Arens
regular }\; \LO Z_r(\pi^{***r***r})={\cal A}^{(4)} \hspace {0.7
cm}(2.1),$$
 and
$$Z_l(\pi^{r***r***})={\cal A}^{(4)}\LO ({\cal A}^{**},\pi^{r***r})
\;\text{is Arens regular}\; \LO Z_r(\pi^{r******r})={\cal
A}^{(4)}\hspace {0.7 cm}(2.2).
$$
On the other hand $\cal A$ is Arens regular if $\cal A^{**}$ is
Arens regular. Let one of the conditions (i), (ii), (iii) or (iv)
holds, then  $\cal A$ is Arens regular; i.e.
$\pi^{r***r}=\pi^{***}.$ Then we have
$$
Z_l(\pi^{r***r***})={\cal A}^{(4)} \LO Z_l(\pi^{******})={\cal
A}^{(4)} \hspace {0.7 cm}(2.3).
$$
(2.1),(2.2) and (2.3) imply that (i),...,(iv) are
equivalent.\hfill$\blacksquare~$

Let ${\cal A}$ be a commutative Banach algebra. Then the following
assertions are equivalent.

(i) ${\cal A} \;\text{is Arens regular}$

(ii) There is $n\in \Bbb N$ such that  ${\cal A}^{(2n)}$ is Arens
regular.

(iii) For every $n\in \Bbb N$ $ {\cal A}^{(2n)}$ is Arens regular.
\paragraph{\bf Theorem 2.2.} Let ${\cal A}$ be a Banach algebra with a bounded right approximate identity, then
$\widehat{{\cal A}^{**}}$ is a left ideal of ${\cal A}^{(4)}$ if and
only if ${\cal A}$ is reflexive.
\paragraph{\bf Proof.} Let $(e_\al)_{\al\in I}$ be a bounded right approximate identity for $\cal A$ with cluster
point $E\in {\cal A}^{**}$. $E$ is a right unit element of ${\cal
A}^{**}$, then $\widehat{E}$ is a right unit element of ${\cal
A}^{(4)}$. If ${\cal A}^{**}$ is a left ideal of ${\cal A}^{(4)}$
then ${\cal A}^{4}={\cal A}^4\widehat{E}=\widehat{{\cal A}^{**}},$
then ${\cal A}$ is reflexive. The converse is
trivial.\hfill$\blacksquare~$

Let $\cal A$ be a Banach algebra in which $\widehat{{\cal A}^{**}}$
be a left ideal of ${\cal A}^{(4)}$. Then we can show that
$\widehat{{\cal A}}$ is a left ideal of ${\cal A}^{**}.$ In the
following we show that the converse of the above statement does not
hold.
\paragraph{\bf Example 2.3.} Let $G$ be a compact topological group, then $L^1(G)$ is an ideal of
$L^1(G)^{**}$. On the other hand $L^1(G)$ is reflexive if and only
if $G$ is finite. Then by above lemma, $L^1(G)^{**}$ is not a left
ideal of $L^1(G)^{(4)}$ when $G$ is infinite.

Let $\pi:{\cal A}\times {\cal A}\longrightarrow {\cal A}$ be the
product of Banach algebra $\cal A$, Dales, Rodriguez and Velasco in
[3] found some necessary and sufficient conditions for Arens
regularity of both $\pi$ and $\pi^{r*}$. They proved the following
theorem that plays a key role in [3]. We will use this theorem to show
that $\pi^{r**}$ is Arens regular if and only if ${\cal A}^{**}$ is
a left ideal of ${\cal A}^{****}$ when  $\pi$ and $\pi^{r*}$ are
Arens regular.
\paragraph{\bf Theorem 2.4.} Let
$f:X\times Y\longrightarrow Z$ be a continuous bilinear map, then
the following assertions are equivalent.

(i) $f$ and $f^*$ are Arens regular.

(ii) $f^{r*r***}=f^{***r*r}.$

(iii) $f^{****}(Z^{***},X^{**})\subseteq \widehat{Y^*}. $
\paragraph{\bf corollary 2.5.} Let $\cal A$ be a Banach algebra
with product
 $\pi:{\cal A}\times {\cal A}\longrightarrow {\cal A}$. Let $\pi$ and $\pi^{r*}$ be Arens
 regular, then $\pi^{r**}$ is Arens regular if and only if ${\cal A}^{**}$ be a left ideal in ${\cal A}^{****}$.
\paragraph{\bf Proof:} Let $f:=\pi^{r*r}$, then $f$ and $f^{r*}$ are Arens
regular if and only if $f^{****}(Z^{***},X^{**})\subseteq \widehat
Y^*, $ we assume that $\pi^{r*r}$ be Arens regular then $\pi^{r**}$
is Arens regular if and only if $\pi^{r*r****}({\cal A}^{****},{\cal
A}^{**})\subseteq \widehat{({\cal A}^{**})}$.\\   For $a''''\in
{\cal A}^{****}$, $a'''\in {\cal A}^{***}$ and $a''\in {\cal
A}^{**}$, we have
\begin{eqnarray*}\langle \pi^{***r*r*}(a'''',a''),a'''\rangle=\langle a'''',
\pi^{***r*}(a''',a'')\rangle\hspace{0,7 cm}(2.4),\end{eqnarray*} and
\begin{eqnarray*}\langle \pi^{******}(a'''',\widehat{a''}),a'''\rangle=\langle a'''',
\pi^{*****}(\widehat{a''},a''')\rangle \hspace{0,7
cm}(2.5).\end{eqnarray*} Also for every $b''\in {\cal A}^{**}$, we
have
\begin{eqnarray*}
\langle \pi^{***r*}(a''',a''),b''\rangle&=&\langle a''',
\pi^{***r}(a'',b'')\rangle =\langle a''', \pi^{***}(b'',a'')\rangle\\
&=&\langle  \pi^{****}(a''', b''), a''\rangle= \langle \widehat{a''}, \pi^{****}(a''', b'')\rangle\\
&=&\langle \pi^{*****}(\widehat{a''},a'''), b''\rangle\hspace{0,7
cm}(2.6).
\end{eqnarray*}
By (2.1), (2.2) and (2.3), we have
$\pi^{r*r****}(a'''',a'')=\pi^{******}(a'''',\widehat{a''}).$ On the
other hand since $\pi$ and $\pi^{r*}$ are Arens regular, then by
theorem 3.1, we have
$\pi^{r*r****}=(\pi^{r*r***})^*=(\pi^{***r*r})^*=\pi^{***r*r*}$ and
since $\pi^{******}$ is the first Arens product of ${\cal
A}^{****}$, then $\pi^{r**}$ is Arens regular if and only if ${\cal
A}^{**}$ is a left ideal in ${\cal A}^{****}$.\hfill$\blacksquare~$

The conditions in above corollary are very strung. If $\cal A$ has
bounded approximate identity, and if $\pi$, $\pi^{r*}$ and
$\pi^{r**}$ are Arens regular, then it is easy to show that $\cal A$
is reflexive. But Arens regularity of $\pi$, $\pi^{r*}$ and
$\pi^{r**}$ dos not implies the rflexivity of $\cal A$, always. For
example if $\cal A$ is a nonreflexive Banach space with trivial
product, then $\pi$, $\pi^{r*}$ and $\pi^{r**}$ are Arens regular,
but $\cal A$ is not reflexive.
\section{Some new ideals}
Let $\cA$ be a Banach algebra. We consider
$$Z_l{(\cal A)}:=\{a\in \cal A~:~
{\cal A}^{**}.\widehat{a}\subseteq \widehat{\cal A} \}\hspace {0.7
cm} (3.1),$$
$$Z_r{(\cal A)}:=\{a\in \cal A~:~
\widehat{a}.{\cal A}^{**}\subseteq \widehat{\cal A} \}\hspace {0.7
cm} (3.2).$$

In this section we show that for a locally compact group $G$, $G$ is
compact if and only if $Z_l(M(G))=L^1(G)$ (if and only if
$Z_r(M(G))=L^1(G)$), and $G$ is finite if and only if
$Z_l(M(G))=M(G)$ (if and only if $Z_r(M(G))=M(G)$). Also we show
that for a semigroup $S$, $sS$ is finite for every $s\in S$ if and
only if $Z_r(l^1(S))=l^1(S)$. So in the case that $l^1(S)$ is
semisimple, we fined the necessary and sufficient condition on $S$
on which $l^1(S)$ to be an annihilator
algebra.\\
It is easy to show that $Z_l(\cal A)$ and $Z_r(\cal A)$ are two
sided ideals of $\cal A$ and $Z_l(\cal A)$ ($Z_r(\cal A)$ ) is a
left (right) ideal of ${\cal A}^{**}.$ Also $Z_l(\cal A)$ ($Z_r(\cal
A)$ ) is the union of all two sided ideals of $\cal A$ which are
left (right) ideals of ${\cal A}^{**}$. Therefore we result the
following.
\paragraph{\bf Theorem 3.1.} Let $\cal A$ be a Banach algebra.
Then the following assertions are equivalent.

(i) $\cal A=Z_l(\cal A)(Z_r(\cal A))$.

(ii) $\cal A$ is a left (right) ideal in ${\cal A}^{**}$.

(iii) For every $a\in \cal A,$ the map $b\longmapsto ab
(b\longmapsto ba):\cal A \longrightarrow \cal A$ is weakly compact.

\paragraph{\bf Proof.} $(i) \LO (ii)$ is straightforward and for $(ii) \LO (iii)$
see Proposition 1.4.13 of [9].\hfill$\blacksquare~$
\paragraph{\bf Theorem 3.2.} Let $\cal A$ be a Banach algebra.
Then
$$Z_l{(\cal A)}=\{a\in \cal A~:~
 \hbox{the mapping } f\longmapsto
af:{\cal A}^{*}\longrightarrow {\cal A}^{*} \hbox{ is~~}
weak^*-weak-\hbox{ continuous }\},$$

and
$$Z_r{(\cal A)}:=\{a\in \cal A~:~
 \hbox{the mapping } f\longmapsto
fa:{\cal A}^{*}\longrightarrow {\cal A}^{*} \hbox{ is~~}
weak^*-weak-\hbox{ continuous }\}.$$
\paragraph{\bf Proof.} We set
$$U=\{a\in \cal A~:~
 \hbox{the mapping } f\longmapsto
af:{\cal A}^{*}\longrightarrow {\cal A}^{*} \hbox{ is~~}
weak^*-weak-\hbox{ continuous }\},$$ and let $a\in U$ and $b\in
{\cal A} $. If
$f_\alpha~\stackrel{weak^*}{-\hspace{-.2cm}-\hspace{-.2cm}\longrightarrow}
f$ in $\cal A^*$, then
$bf_\alpha~\stackrel{weak^*}{-\hspace{-.2cm}-\hspace{-.2cm}\longrightarrow}
bf$ in $\cal A^*$.  By definition of $U,$
$abf_\alpha~\stackrel{weak}{-\hspace{-.2cm}-\hspace{-.2cm}\longrightarrow}
abf$ in $\cal A^*$. Thus $ab\in U$. On the other hand since $a\in
U$, then
$af_\alpha~\stackrel{weak}{-\hspace{-.2cm}-\hspace{-.2cm}\longrightarrow}
af$ in $\cal A^*$ and because $\cal A^*$ is a dual Banach space,
then every bounded linear map on $\cal A^*$ is weak-weak continuous,
therefore we have
$baf_\alpha~\stackrel{weak}{-\hspace{-.2cm}-\hspace{-.2cm}\longrightarrow}
baf$ in $\cal A^*$. Thus $ba \in U$ and U is a two sided deal in
$\cal A$. Let now $a''\in {\cal A}^{**}$, $a\in U$ and
$f_\alpha~\stackrel{weak^*}{-\hspace{-.2cm}-\hspace{-.2cm}\longrightarrow}
f$ in $\cal A^*$. Then
$af_\alpha~\stackrel{weak}{-\hspace{-.2cm}-\hspace{-.2cm}\longrightarrow}
af$ in $\cal A^*$ and we have $$ \lim_\alpha\langle a''\widehat{a} ,
f_\alpha \rangle=\lim_\alpha\langle a'',af_\alpha\rangle =\langle
a'', af \rangle =\langle a''\widehat{a}, f \rangle.
$$
Thus $a''\widehat{a}:{\cal A}^{**}\longrightarrow \Bbb C$ is
$weak^*$-$weak^*-$ continuous then $a''\widehat{a} \in {\widehat
{\cal A}}$
and $U\subseteq Z_l(\cal A)$.\\
Let now $a''\in {\cal A}^{**}$, $a\in Z_l(\cal A)$ then
$a''\widehat{a} \in {\widehat {\cal A}}$. Suppose
$f_\alpha~\stackrel{weak^*}{-\hspace{-.2cm}-\hspace{-.2cm}\longrightarrow}
f$ in $\cal A^*$. Then
\begin{eqnarray*}
\lim_\alpha\langle a'' , af_\alpha \rangle=\lim_\alpha\langle
a''\widehat{a},f_\alpha\rangle =\langle a''\widehat{a}, f \rangle
=\langle a'', af \rangle.
\end{eqnarray*}
Therefore
$af_\alpha~\stackrel{weak}{-\hspace{-.2cm}-\hspace{-.2cm}\longrightarrow}
af$ in $\cal A^*$ and $Z_l(\cal A)\subseteq U$. Similarly we can
prove the argument for $Z_l(\cal A)$.\hfill$\blacksquare~$
\paragraph{\bf  Example 3.3.} Let $\cal
A=l^1(\Bbb N)$ with product $fg=f(1)g~~~(f,g\in l^1(\Bbb N))$. Then
$\cal A$ is a Banach algebra with $l^1$-norm. $\cal A$ is a left
ideal of ${\cal A}^{**}$ and we have  $\cal A=Z_l(\cal A)(Z_r(\cal
A))$. On the other hand $Z_r(\cal A)=\{f\in \cal A: f(1)=0\}$. Thus
$Z_l(\cal A)$ and $Z_r(\cal A)$ are different.

Let $\cal A$ be a Banach algebra and let $\cal A^{**}$ be a left
(right) ideal of $\cal A^{(4)},$ then $\cal A$ is a left (right)
ideal of $\cal A^{**}.$ Therefore we have the following .

(i) If $Z_l(\cal A^{**})=\cal A^{**},$ then $Z_l(\cal A)=\cal A.$

(ii) If $Z_r(\cal A^{**})=\cal A^{**},$ then $Z_r(\cal A)=\cal A.$
\paragraph{\bf Theorem 3.4.} Let $G$ be a locally compact
group. Then $Z_l(L^1(G))=Z_l(M(G))$ and $Z_r(L^1(G))=Z_r(M(G)).$
\paragraph{\bf  Proof.} Because  $L^1(G)$  has a bounded approximate
identity, then by Cohen factorization theorem, $(L^1(G))^2=L^1(G).$
On the other hand $L^1(G)$ is a two sided ideal of $M(G)$. Therefore
every ideal of $L^1(G)$ is an ideal of $M(G)$. Let
$\pi:L^1(G)\longrightarrow M(G)$ be the inclusion map, then
$\pi^{''}(l^1(G))^{**}$ is a two sided ideal of $M(G)^{**}$. Thus
every left (right or two sided) ideal of $\pi^{''}(l^1(G))^{**}$ is
a left (right or two sided) ideal of $M(G)^{**}$. Then
$Z_l(L^1(G))\subseteq Z_l(M(G)).$ We have to show that
$Z_l(M(G))\subseteq Z_l(L^1(G)).$ To this end let $(e_\alpha)$ be a
bounded approximate identity of $L^1(G)$ with bound 1 and  with a
cluster point $E\in {L^1(G)}^{**}.$ Then the mapping $m\longmapsto
({\pi}^{''}(E))\widehat{m}:M(G)\longrightarrow
{\pi}^{''}(L^1(G))^{**}$ is isometric embedding. We denote this map
with $\Gamma_E$. Since the restriction of $\Gamma_E$ to $L^1(G)$ is
identity map, then $\Gamma_E(m) \in \widehat {{\pi}(L^1(G))}$ if and
only if $m\in L^1(G)$. Let now $m\in Z_l(M(G))$ then
$M(G)^{**}\widehat{m}\subseteq \widehat {M(G)}$. Thus
${\pi}^{''}(L^1(G))^{**}\widehat{m}\subseteq \widehat {M(G)}$. Since
${\pi}^{''}(L^1(G))^{**}$ is an ideal of $M(G)^{**}$, we have
${\pi}^{''}(L^1(G))^{**}\widehat{m}\subseteq [(\widehat {M(G)})\cap
{\pi}^{''}(L^1(G))^{**}]=\widehat {L^1(G)}$ (see corollary 3.4 of
[6]). Therefore $\Gamma_E(m) \in {\pi}^{''}(L^1(G))^{**}$ and $m\in
L^1(G)$. This conclude that $Z_l(M(G))\subseteq Z_l(L^1(G)).$
Similarly we can show that $Z_r(M(G))=
Z_r(L^1(G)).$\hfill$\blacksquare~$
\paragraph{\bf  Corollary 3.5.} For a locally compact group $G$
the following assertions are equivalent.

(i) $G$ is compact.

(ii) $Z_l(L^1(G))=L^1(G)$ ($~~Z_r(L^1(G))=L^1(G)$).

(iii) $Z_l(M(G))=L^1(G)$ ($~~Z_r(M(G))=L^1(G)).$

\paragraph{\bf Theorem 3.6.} Let $S$ be a semigroup, then

(i) $sS$ is finite for every $s\in S$.

(ii) $Z_r(l^1(S))=l^1(S)$.
\paragraph{\bf Proof.} (i)$\Lo$(ii). We have to show that for every $a\in l^1(S)$, $\la_a:l^1(S)\lo
l^1(S)$ is compact operator. To this end, let $a\in l^1(S)$, then
$a=\sum_{n=1}^\infty a_ns_n$ when $c_n=a(s_n)$. Since $s_nS$ is
finite for every $n$, then $\la_{s_n}(l^1(S))$ is a finite dimension
subspace of $l^1(S)$. Thus $\la_{s_n}$ is compact operator on
$l^1(S)$. Therefore the operator $\sum_{n=1}^k c_n\la_{s_n}$ is
compact for every $k\in\BN$. But $\la_a=\sum_{n=1}^\infty
a_n\la_{s_n}$, then by VI. 5.3. of [5], the set of compact operators
is closed in the uniform operator topology of $BL(X,Y)$ and we get
$\la_a$ is a compact operator on $l^1(S)$.

(ii)$\Lo$(i). Let $s_0\in S$ and $s_0S$ be infinite. Then there
exists $\{u_n\}_{n\in\BN}\subs S$ such that $s_0u_n\neq s_0u_m$ when
 $n\neq m$. Then $\la_{s_0}$ is isometric on an infinite dimension
subspace of $l^1(s)$. i.e. $\la_{s_0}$ is not
compact.\hfill$\blacksquare~$
\paragraph{\bf Corollary 3.7.} Let $l^1(S)$ be semisimple. Then the following assertions are equivalent.

(i) $Z_r(l^1(S))=Z_l(l^1(S))=l^1(S)$ and $S=\{st: s,t\in S\}.$

(ii) $l^1(S)$ is an annihilator algebra.
\paragraph{\bf Proof.} (i)$\Lo$(ii). Let $s\in S$ and let (i) holds. Then $SsS$ is
finite therefore $l^1(S) sl^1(S)$ is finite dimension. Since
$l^1(S)$ is semisimple and $l^1(S) sl^1(S)$ is an ideal of $l^1(S),$
then $l^1(S) sl^1(S)$ is semisimple finite dimension ideal of
$l^1(S).$ Therefore $l^1(S) sl^1(S)$ is isomorphic with the direct
sum of full matrix algebras. Now, let $P\in S.$ Then
 $P=s_1s_2$ for some $s_1$ and $s_2$ in $S$ and we have
$s_2=t_1t_2$ for $t_1,t_2\in S$. Thus $P=s_1t_1t_2\in St_1S$ for
some $t_1\in S$. On the other hand for each $a\in l^1(S)$ we have
$a=\sum_{n=1}^\infty C_nP_n$ where $P_n\in S$ and $a(P_n)=C_n$. We
get that $l^1(S)$ is the topological  sum of full matrix algebras,
and by 2.8.29 of [8], $l^1(S)$ is an annihilator algebra.
(ii)$\Lo$(i). Since $l^1(S)$ is semisimple annihilator algebra, then
by theorem 3.1. of [10] $\widehat{l^1(S)}$ is a two sided ideal of
$(l^1(S))^{**}$. Then by above theorem, $Ss$ and $sS$ are finite for
every $s\in S$. To prove $S=\{st: s,t\in S\}$ we have
$l^1(S)l^1(S)\subs l^1(S^2)$ where $l^1(S^2)$ is a closed two sided
ideal of $l^1(S)$. Now $l^1(S)$ is an annihilator  algebra, then
\begin{eqnarray*} ran(l^1(S))=\{0\} &\Lo& ran(l(S))^2=\{0\}\\
&\Lo&ran(l^1(S^2))=\{0\} \\
&\Lo& l^1(S^2)=l^1(S)\\
&\Lo&S^2=S.\end{eqnarray*}\hfill$\blacksquare~$


\end{document}